\documentclass[12pt]{article}
\usepackage{stmaryrd}
\usepackage{mathrsfs}
\usepackage{amsmath}
\usepackage{amsmath,amsthm,amssymb}
\usepackage{amsfonts}
\usepackage{latexsym}
\date{}

\begin{document}
\title{Representation dimension of $m$-replicated algebras$^\star$}
\author{{\small  Hongbo Lv, Shunhua Zhang$^{*}$}\\
{\small  Department of Mathematics,\ Shandong University,\ Jinan
250100, P. R. China }\\
{\small Dedicated to Professor Shaoxue Liu on the occasion of his
eightieth birthday}}

\pagenumbering{arabic}

\maketitle
\begin{center}
 \begin{minipage}{120mm}
   \small\rm
   {\bf  Abstract}\ \ Let $A$ be a finite dimensional hereditary algebra over an algebraically
closed field and $A^{(m)}$ the $m$-replicated algebra of  $A$. We
prove that the representation dimension of $A^{(m)}$ is at most
three, and that the dominant dimension of $A^{(m)}$ is at least $m$.

\vskip 0.1in

 {\bf Key words: $m$-replicated algebras, representation dimension, dominant dimension}

\end{minipage}
\end{center}
\footnote {MSC(2000): 16E10, 16G10}

\footnote{ $^\star$Supported by the NSF of China (Grant No.
10771112).}

\footnote{ $^{*}$Corresponding author.}

\footnote{ {\it Email addresses}: lvhongbo356@163.com(H.Lv), \
shzhang@sdu.edu.cn(S.Zhang)}

\section {Introduction}

\vskip 0.2in

Representation dimension of Artin algebras was introduced by M.
Auslander in [A] to study the connection of arbitrary Artin algebras
with representation finite Artin algebras. It gives a reasonable way
to understand how far an Artin algebra is from being representation
finite type by measuring the global dimension of all endomorphism
rings of modules which are both generators and cogenerators. M.
Auslander proved that an Artin algebra $\Lambda$ is representation
finite if and only if its representation dimension is at most two.

\vskip 0.2in

An interesting relationship between the representation dimension and
the finitistic dimension conjecture was shown by K. Igusa and G.
Todorov [IT], which is, if the representation dimension of an
algebra is at most three, then its finitistic dimension is finite.
Furthermore, in [ZZ], A.Zhang and S.Zhang proved that if
quasi-hereditary algebras have representation dimensions at most
three, then the finitistic dimension conjecture holds. Recently,
many important classes of algebras, such as tilted algebras, laura
algebras, trivial extensions of iterated tilted algebras etc., have
been shown to have representation dimensions at most three, see
[APT, CP, EHIS, X].

\vskip 0.2in

Let $A$ be a hereditary algebra over an algebraically closed field
$k$, and $A^{(m)}$ be the $m$-replicated algebra of $A$.  This kinds
of algebras give a one-to-one correspondence between the basic
tilting $A^{(m)}$-modules and the basic tilting objects in
$m$-cluster category $\mathscr{C}_m(A)$, see [ABST1, ABST2].  This
motivates further investigation on tilting modules about this kinds
of algebras, see [LLZ,Z1,Z2] for details. In this paper, we prove
that the representation dimension of $A^{(m)}$ is at most three, and
that the dominant dimension of  $A^{(m)}$ is at least $m$.

\vskip 0.2in

The following theorems is the main results of this paper.

\vskip 0.2in

{\bf Theorem 1.}\ {\it Let $A^{(m)}$ be the $m-$replicated algebra
of a hereditary algebra $A$. Then its representation dimension is at
most three.}

\vskip 0.2in

{\bf Theorem 2.}\ {\it Let $A^{(m)}$ be the $m-$replicated algebra
of a hereditary algebra $A$. Then its dominant dimension is at least
$m$.}

\vskip 0.2in

This paper is arranged as follows. In section 2, we fix the
notations and collect necessary definitions and basic facts. Section
3 is devoted to the proof of our main results.

\vskip 0.2in

\section { Preliminaries}

\vskip 0.2in

Let $\Lambda$ be a basic connected algebra over an algebraically
closed field $k$. We denote by mod $\Lambda$ the category of all
finitely generated right $\Lambda$-modules and by ind $\Lambda$ a
full subcategory of mod $\Lambda$ containing exactly one
representative of each isomorphism class of indecomposable
$\Lambda$-modules. The  bounded derived category of mod $\Lambda$ is
denoted by $\mathcal{D}^{b}({\rm mod}\ \Lambda)$ and the shift
functor by [1]. For a $\Lambda$-module $M$, we denote by add $M$ the
full subcategory of $\rm{mod}\ \Lambda$ whose objects are the direct
summands of finite direct sums of copies of $M$ and by
$\Omega_{\Lambda}^{-1}M$ the first cosyzygy which is the cokernel of
an injective envelope $M\hookrightarrow I.$ The projective dimension
of $M$ is denoted by {pd $M$}, the global dimension of $\Lambda$ by
gl.dim $\Lambda$ and the Auslander-Reiten translation of $\Lambda$
by $\tau_\Lambda$. $D={\rm Hom}_k(-,\ k)$ is the standard duality
between mod $\Lambda$ and mod $\Lambda^{op}$. For further
definitions and facts about mod $\Lambda$, we refer to [ARS, Rin].

\vskip 0.2in

Let $\mathcal{C}$ be a full subcategory of mod $\Lambda$,
$C_{M}\in\mathcal{C}$ and $\varphi :C_M\longrightarrow M$ with
$M\in$ mod $\Lambda$. The morphism $\varphi$ is a right
$\mathcal{C}$-approximation of $M$ if the induced  morphism ${\rm
Hom}(C,C_{M})\longrightarrow {\rm Hom}(C,M)$ is surjective for any
$C\in\mathcal{C}$. A minimal right $\mathcal{C}$-approximation of
$M$ is a right $\mathcal{C}$-approximation which is also a right
minimal morphism, i.e., its restriction to any nonzero summand is
nonzero. The subcategory $\mathcal{C}$ is called contravariantly
finite if any module $M\in$ mod $\Lambda $ admits a (minimal) right
$\mathcal{C}$-approximation. If $\mathcal{C}$ is also closed under
extensions and $K$ is the kernel of a minimal right
$\mathcal{C}$-approximation of $M$, then Wakamatsu's Lemma states
that ${\rm Ext}^{1}_{\mathcal{C}}(L, K) = 0$ for all $L\in
\mathcal{C}$ (see [W]). The notions of (minimal) left
$\mathcal{C}$-approximation and of covariantly finite subcategory
are dually defined. It is well known that add $M$ is both a
contravariantly finite subcategory and a covariantly finite
subcategory.

\vskip 0.2in

Let $M$ and $N$ be two indecomposable $\Lambda$-modules. A path from
$M$ to $N$ in ind $\Lambda$ is a sequence of non-zero morphisms
$$
M=M_0\stackrel{f_1} \longrightarrow M_1\stackrel{f_2} \longrightarrow\cdots
  \stackrel{f_t} \longrightarrow M_t =N
$$
with all $M_i$ in  ind $\Lambda$. Following [Rin], we denote the
existence of such a path by $M\leq N$. We say that $M$ is a
predecessor of $N$ (or that $N$ is a successor of $M$). More
generally, if $S_1$ and $S_2$ are two sets of modules, we
 write $S_1\leq S_2$ if every module in $S_2$ has a predecessor in
 $S_1$, every module in $S_1$ has a successor in
 $S_2$, no  module in $S_2$ has a successor in
 $S_1$ and no  module in $S_1$ has a predecessor in
 $S_2$. The notation $S_1<S_2$ stands for $S_1\leq S_2$ and $S_1\cap S_2 \neq \emptyset.$

\vskip 0.2in

Let $\hat{\Lambda }$ be the repetitive algebra of $\Lambda$. Then
$\hat{\Lambda }$ is the infinite matrix algebra
$$
\hat{\Lambda} =\begin{pmatrix}
                                \ddots &  & 0 &  &  \\
                                 & \Lambda_{i-1} &  &  &  \\
                                & Q_{i} & \Lambda_{i} &  &  \\
                                &  & Q_{i+1} & \Lambda_{i+1} &  \\
                                & 0 & & \ddots \\
                             \end{pmatrix}
$$
where matrices have only finitely many non-zero coefficients,
$\Lambda_{i}=\Lambda$ and $Q_{i}= {_{\Lambda}D\Lambda_{\Lambda}}$
for all $i\in \mathbb{Z}$, all the remaining coefficients are zero
and multiplication is induced from the canonical isomorphisms
$\Lambda \otimes_\Lambda D\Lambda\cong\!
_{\Lambda}D\Lambda_{\Lambda}\cong D\Lambda\otimes_\Lambda \Lambda $
and the zero morphism $D\Lambda\otimes_\Lambda D\Lambda
\longrightarrow 0$ (see [HW]). From [H], we have the following
important connection between $\hat{\Lambda} $ and
$\mathcal{D}^{b}({\rm mod}\ \Lambda)$.

\vskip 0.2in

{\bf  Lemma 2.1.}\ {\it Let $\Lambda$ be of finite global dimension,
then the derived category $\mathcal{D}^{b}({\rm mod}\ \Lambda)$ is
equivalent, as a triangulated category, to the stable module
category ${\rm \underline{mod}}\ \hat{\Lambda}$.}

\vskip 0.2in

The right repetitive algebra $\Lambda'$ of $\Lambda$, introduced in
[ABM], is the quotient of $\hat{\Lambda}$ defined as follows:
$$\Lambda'=\begin{pmatrix}
             \Lambda_{0} &  &  & 0 \\
             Q_{1} & \Lambda_{1} &  &  \\
              & Q_{2} & \Lambda_{2} &  \\
             0 &  & \ddots & \ddots \\
           \end{pmatrix}
$$
where, as above, $\Lambda_{i}=\Lambda$ and
$Q_{i}={_{\Lambda}D\Lambda_{\Lambda}}$ for all $i\geq 0$. We denote
by $\Sigma_0$ the set of all non-isomorphic indecomposable
 projective $\Lambda$-modules and set
 $\Sigma_k=\Omega_{\Lambda'}^{-k}\,\Sigma_0$ for $k\geq 0$.

\vskip 0.2in

Then, by [AI], the $m-$replicated algebra $\Lambda^{(m)}$ of
$\Lambda$ is defined as the quotient of the right repetitive algebra
$\Lambda'$, hence of the repetitive algebra $\hat{\Lambda}$, which
is,
$$\Lambda^{(m)}=\begin{pmatrix}
                  \Lambda_{0} &  &  &  & 0 & \\
                  Q_{1} & \Lambda_{1} & &  &  &  \\
                   & Q_{2} & \Lambda_{2} &  &  &  \\
                 &  &  & \ddots & \ddots &  \\
                  &  0&  &  & Q_{m} & \Lambda_{m} \\
                \end{pmatrix}.
$$

\vskip 0.2in

If $m=1$, then $\Lambda^{(1)}$ is the duplicated algebra of
$\Lambda$ (see [ABST1]). Also from [AI], we have that
$$
m+{\rm gl.dim }\ \Lambda \leq {\rm gl.dim }\ \Lambda^{(m)}\leq(m+1){\rm gl.dim }\ \Lambda+m.
$$

\vskip 0.2in

  The next lemma is from [ABST2].

\vskip 0.2in

{\bf Lemma 2.2.}\ {\it Let $A$ be hereditary. Then

\vskip 0.1in

 {\rm (1)}The standard embeddings ${\rm ind}\
A_{i}\hookrightarrow {\rm ind }\ A^{(m)}$ (where $0\leq i\leq m$)
and ${\rm ind }\ A^{(m)} \hookrightarrow {\rm ind }\ A'$ are full,
exact, preserve indecomposable modules, almost split sequences and
irreducible morphisms.

\vskip 0.1in

{\rm(2)}\ Let $M$ be an indecomposable $A'$-module which is not
projective and $k\geq 1$. Then the followings are equivalent:

{\rm(a)}\ {\rm pd} $M =k$,

{\rm(b)}\ there exists $N \in $ {\rm ind}\ $A$ such that $M\cong
\tau^{-1}_{A'} \Omega^{-(k-1)}_{A'}N$,

{\rm(c)}\ $\Sigma_{k-1}< M\leq \Sigma_{k}.$

\vskip 0.1in

{\rm(3)}\ Let $M$ be an indecomposable $A^{(m)}$-module which is not
in ${\rm ind}\ A= {\rm ind}\ A_0$. Then its projective cover in {\rm
mod} $A^{(m)}$ is projective-injective and coincides with its
projective cover in {\rm mod} $A'$.

\vskip 0.1in

{\rm(4)}\ Let $M$ be an $A^{(m)}$-module having all
projective-injective indecomposable modules as direct summands. For
$A^{(m)}$-module $X$, if $X$ has a projective-injective projective
cover, then a minimal right {\rm add} $M$-approximation of $X$ is an
epimorphism.}

\vskip 0.2in

We refer to [A] for the original definition, and we would rather use
the following characterisation, which was proved to be equivalent to
the original one in [A].

\vskip 0.2in

{\bf Definition 2.3.}\ Let $\Lambda$ be a non-semisimple Artin
algebra. The representation dimension rep.dim $\Lambda$ of $\Lambda$
is the infimum of the global dimensions of the algebras End $M$,
where $M$ is a generator and a cogenerator of mod $\Lambda$.

\vskip 0.2in

The next lemma in [A, CP, EHIS, X] is well known.

\vskip 0.2in

{\bf Lemma 2.4.}\ {\it Let $\Lambda$ be a non-semisimple Artin
algebra, $n$ be a positive integer, and $M$ be a
generator-cogenerator of {\rm mod} $\Lambda$. Then {\rm gl.dim End}
$M \leq n+1$ if and only if for each $\Lambda$-module $X$, there
exists an exact sequence
$$0\rightarrow M_{n}\rightarrow \cdots \rightarrow M_1\rightarrow X\rightarrow 0$$
with $M_{i}$ in {\rm add} $M$ for all $i$, such that the induced
sequence
$$0\rightarrow{\rm Hom}_{\Lambda}(L,M_{n})\rightarrow \cdots \rightarrow{\rm Hom}_{\Lambda}(L,M_{1})
\rightarrow{\rm Hom}_{\Lambda}(L,X)\rightarrow0$$ is exact for all
$L$ in {\rm add} $M.$ In particular, {\rm rep.dim} $\Lambda \leq
n+1$. }

\vskip 0.2in

We refer to [ARS] for the following definition.

\vskip 0.2in

{\bf Definition 2.5.}\ Let $\Lambda$ be a non-semisimple Artin
algebra. The dominant dimension dom.dim $\Lambda$ of $\Lambda$ is
the supremum of all $n\in \mathbb{N}$ having the property that if
$$
0\rightarrow \Lambda \rightarrow I_1\rightarrow I_2 \rightarrow
\cdots \rightarrow I_n\rightarrow \cdots
$$
is the minimal injective resolution of $\Lambda$,  the $T_j$ is
projective for all $j\leq n$.

\vskip 0.2in

\section {The representation dimension of $A^{(m)}$}

\vskip 0.2in

  In this section, we will prove our main result.

\vskip 0.1in

Let $\Lambda$ be a, not necessarily hereditary, finite dimensional
algebra and $\hat{\Lambda}$ the repetitive algebra of $\Lambda$. For
$\hat{\Lambda}-$modules $X, Y$, we define $$\varphi: {\rm
Hom}_{\hat{\Lambda}}(Y, \Omega^{-1}_{\hat{\Lambda}}X)
\longrightarrow{\rm Ext}^{1}_{\hat{\Lambda}}(Y, X)$$
 as follows: for
$h\in{\rm Hom}_{\hat{\Lambda}}(Y, \Omega^{-1}_{\hat{\Lambda}}X)$, we
have the following commutative diagram with exact rows
$$\begin{array}{ccccccccc}
  0 & \longrightarrow & X &\longrightarrow
  & E & \longrightarrow &Y&  \longrightarrow & 0 \\
    &  &\parallel &  & \Big\downarrow
&  & \Big\downarrow \vcenter{%
\rlap{$\scriptstyle{h}$}}&  &  \\
  0 & \longrightarrow & X& \longrightarrow  & I_{X}
  & \longrightarrow & \Omega^{-1}_{\hat{\Lambda}}X & \longrightarrow & 0, \\
\end{array}$$where the right square is a pullback and $I_{X}$
is the injective envelope of $X$ and also the projective cover of
$\Omega^{-1}_{\hat{\Lambda}}X$. Define $\varphi(h)$ as the short
exact sequence
$$0\longrightarrow X \longrightarrow E\longrightarrow Y\longrightarrow0.$$
By the uniqueness of pullback, $\varphi$ is well defined and an
epimorphism. The next lemma shows that $\varphi$ can induce an
isomorphism between ${\underline{\rm Hom}}_{\hat{\Lambda}}(Y,
\Omega^{-1}_{\hat{\Lambda}}X)$ and ${\rm Ext}^{1}_{\hat{\Lambda}}(Y,
X)$.

\vskip 0.2in

{\bf Lemma 3.1.}\ {\it Let $h\in {\rm Hom}_{\hat{\Lambda}}(Y,
\Omega^{-1}_{\hat{\Lambda}}X)$. Then
   $h$  factors through projective-injective
$\hat{\Lambda}$-modules if and only if $\varphi(h)$ is zero in ${\rm
Ext}^{1}_{\hat{\Lambda}}(Y, X)$.}

\vskip 0.1in

{\bf Proof.}\ If $\varphi(h)$ is zero in ${\rm
Ext}^{1}_{\hat{\Lambda}}(Y, X)$, we can write $\varphi(h)$ as $$0
\longrightarrow X\longrightarrow X\oplus Y
\stackrel{f}\longrightarrow Y\longrightarrow 0,$$ where $f$ is a
retraction. Then the following commutative diagram with exact rows
$$\begin{array}{ccccccccc}
  0 & \longrightarrow & X &\longrightarrow
  & X\oplus Y & \underset{f'}{\stackrel{f}\rightleftharpoons} &Y&  \longrightarrow & 0 \\
    &  &\parallel &  & \Big\downarrow \vcenter{%
\rlap{$\scriptstyle{g}$}}
&  & \Big\downarrow \vcenter{%
\rlap{$\scriptstyle{h}$}}&  &  \\
  0 & \longrightarrow & X& \longrightarrow  & I_{X}
  & \stackrel{\pi}\longrightarrow & \Omega^{-1}_{\hat{\Lambda}}X & \longrightarrow & 0 \\
\end{array}$$
yields that $h=\pi gf'$, which implies that $h$ can factor through a
projective-injective $\hat{\Lambda}-$module.

Conversely, assume that $h$ can factor through a
projective-injective $\hat{\Lambda}-$module. Let $\varphi(h)$ be the
short exact sequence
$$
0\longrightarrow X\longrightarrow E \longrightarrow
Y\longrightarrow 0.
$$
Therefore we have, by [H], a triangle in the
stable module category $\underline{{\rm mod}}\ \hat{\Lambda}$
$$X\longrightarrow E \longrightarrow Y \stackrel{\underline{h}}\longrightarrow
\Omega^{-1}_{\hat{\Lambda}}X,$$ where $\underline{h}$ is zero in
${\rm \underline{Hom}}_{\hat{\Lambda}}(Y,
\Omega^{-1}_{\hat{\Lambda}}X)$, by our assumption. Therefore this
triangle is split, and thus $X$ and $Y$ are both direct summands of
$E$, which implies $\varphi(h)$ is zero. This completes the proof.
\hfill$\Box$

\vskip 0.2in

Assume from now on that $A$ is a hereditary algebra, and that $A'$
is the right repetitive algebra of $A$ and $A^{(m)}$ is the
$m$-replicated algebra.

\vskip 0.2in

{\bf Lemma 3.2.}\ {\it Let $X\in {\rm ind}\ A$, and $\alpha \in {\rm
Hom}_{\hat{A}}(\Omega^{-i}_{\hat{A}}A, \Omega^{-j}_{\hat{A}}X).$ If
$i<j$, then $\alpha$ factors through a projective-injective
$\hat{A}$-module. }

\vskip 0.1in

{\bf Proof.}\ The statement follows easily from that
$$\begin{array}{l}
  {\rm\underline{Hom} }_{\hat{A}}(\Omega^{-i}_{\hat{A}}A,
\Omega^{-j}_{\hat{A}}X) \\
  \cong {\rm Hom}_{\mathcal{D}^{b}({\rm mod}A
)}(A[i], X[j]) \\
  ={\rm Hom}_{\mathcal{D}^{b}({\rm mod}A )}(A,
X[j-i]) \\
  ={\rm Ext}_{A}^{j-i}(A, X) \\
  =0 .
\end{array}$$
 \hfill$\Box$

\vskip 0.2in

Now let gl.dim $A^{(m)}=t.$ Recall that $m+1\leq t\leq 2m+1$,
$\Sigma_0$ is the set of all non-isomorphic indecomposable
 projective $A$-modules and
that $\Sigma_{k}=\Omega_{A'}^{-k}\,\Sigma_0$. By Lemma 2.4 and the
Auslander-Reiten quiver of $A^{(m)}$, $\Sigma_{t-1}\bigcap {\rm ind
}\ A_m \neq \emptyset$ or $\tau^{-1}_{A^{(m)}}\Sigma_{t-1}\bigcap
{\rm ind }\ A_m \neq \emptyset.$ We denote by $U_{k}$ the direct sum
of all the indecomposable modules in $\Sigma_{k}\bigcap {\rm ind}\
A^{(m)}$ and by $P$ the direct sum of all the projective-injective
indecomposable $A^{(m)}$-modules.

\vskip 0.2in

{\bf Theorem 3.3.}\ {\it Let $M=A\oplus DA_{m}\oplus P\oplus
\underset{k=1}{\overset{t-1}{\oplus}}U_k$. Then {\rm gl.dim} ${\rm
End}_{A^{(m)}}(M)\leq 3$. In particular, {\rm rep.dim} $A^{(m)}\leq
3$.}

\vskip 0.1in

{\bf Proof.}\ By Lemma 2.4, it suffices to find, for each
indecomposable $A^{(m)}$-module $X$, a short exact sequence
$$0\longrightarrow M_2\longrightarrow M_1\longrightarrow X\longrightarrow 0$$
with $M_1, M_2 \in {\rm add}\ M$, such that the induced sequence
$$0\longrightarrow {\rm Hom}_{A^{(m)}}(L, M_2)\longrightarrow {\rm Hom}_{A^{(m)}}(L, M_1)
\longrightarrow {\rm Hom}_{A^{(m)}}(L, X)\longrightarrow 0$$ is
exact for all $L\in {\rm add}\ M$.

\vskip 0.1in

It is clear that we can assume that $X$ is not in add $M$. Suppose
firstly that $\Sigma_0<X<\Sigma_1.$ Consider the minimal projective
resolution
$$0 \longrightarrow P_2\longrightarrow P_1 \stackrel{f}\longrightarrow X\longrightarrow 0. $$
Since $X\in {\rm ind}\ A$ and ${\rm Hom}_{A^{(m)}}(DA_{m}\oplus
P\oplus \underset{k=1}{\overset{t-1}{\oplus}}U_k, X)=0$, it is easy
to see that $f$ is a right add $M$-approximation of $X$. Therefore
the induced sequence
$$0\longrightarrow {\rm Hom}_{A^{(m)}}(L, P_2)\longrightarrow {\rm Hom}_{A^{(m)}}(L, P_1)
\longrightarrow {\rm Hom}_{A^{(m)}}(L, X)\longrightarrow 0$$ is
exact for all $L\in {\rm add}\ M$.

\vskip 0.1in

Assume now $\Sigma_i<X<\Sigma_{i+1}$, $1\leq i\leq t-2.$ Then $X$ is
not in ind $A$ and by Lemma 2.2(2), $X$ has a projective-injective
projective cover. It follows from Lemma 2.2(4) that a minimal right
add $(U_i\oplus P)$-approximation of $X$ is an epimorphism. Consider
the following short exact sequence
$$0 \longrightarrow K \longrightarrow M_1 \stackrel{g}\longrightarrow X\longrightarrow 0,\ \ \ \ \ \ (\ast)$$
where $g$ is a minimal right add $(U_i\oplus P)$-approximation of
$X$ and $K$ is the kernel of $g$. By Lemma 2.2(2), there is an
indecomposable $A$-module $Y$ such that
$X\cong\Omega^{-i}_{A^{(m)}}Y$. Then $g$ is also a minimal right add
$M$-approximation of $X$ by Lemma 3.2 and the fact that ${\rm
Hom}_{A^{(m)}}(DA_{m}\oplus\underset{i<j\leq t-1}{\oplus}U_j, X)=0.$
Since the short exact sequence $(\ast)$ is not split, $K$ is not
projective-injective. Clearly, $K\leq\Sigma_i.$ Let $N$ be a
non-projective-injective indecomposable direct summand of $K$ and
assume that $\Sigma_l<N<\Sigma_{l+1}$, for some $0\leq l\leq i-1$.
Then there exists, by Lemma 2.2(2), an indecomposable $A$-module
$N'$ such that $N\cong\Omega^{-l}_{A^{(m)}}N'$. It follows from
Wakamatsu's Lemma that ${\rm
Ext}_{A^{(m)}}^{1}(\Omega^{-(l+1)}_{A^{(m)}}A, K)=0$. Then we have
that $$\begin{array}{rlr}
   0& ={\rm
Ext}_{A^{(m)}}^{1}(\Omega^{-(l+1)}_{A^{(m)}}A, K) \\
   & \cong{\rm
Ext}_{\hat{A}}^{1}(\Omega^{-(l+1)}_{\hat{A}}A, K)&\ \ \ \ \ \ (1) \\
   & \cong{\underline{\rm
Hom}}_{\hat{A}}(\Omega^{-(l+1)}_{\hat{A}}A, \Omega_{\hat{A}}^{-1}K)&\ \ \ \ \ \ (2) \\
   & \cong{\rm
Hom}_{\mathcal{D}^{b}({\rm mod}A)}(A[l+1], K[1]) &\ \ \ \ \ \ (3)\\
   & \cong{\rm
Hom}_{\mathcal{D}^{b}({\rm mod}A)}(A, K[-l])&\ \ \ \ \ \ (4)
\end{array}$$
where (1) follows from Lemma 2.2(1), (2) is from Lemma 3.1, (3)
follows from Lemma 2.1 and (4) holds because [1] is a
selfequivalence of $\mathcal{D}^{b}({\rm mod}\ A).$ In particular,
$${\rm Hom}_{\mathcal{D}^{b}({\rm mod}A)}(A, N[-l])\cong{\rm
Hom}_{A}(A, N')=0,$$ and hence $N'=N=0,$ which implies that $K\in
{\rm add}(\underset{0\leq j\leq i}{\oplus}U_j\oplus P)$ and hence
$K\in$ add $M$. Since $g$ is a right add $M$-approximation of $X$,
the induced sequence
$$0\longrightarrow {\rm Hom}_{A^{(m)}}(L, K)\longrightarrow {\rm Hom}_{A^{(m)}}(L, M_1)
\longrightarrow {\rm Hom}_{A^{(m)}}(L, X)\longrightarrow 0$$ is
exact for all $L\in {\rm add}\ M$.

\vskip 0.1in

Finally, if $\Sigma_{t-1}< X<DA_m$, by the same argument as above,
we consider the short exact sequence $$0\longrightarrow
K'\longrightarrow M_{1}^{'}\stackrel{h}\longrightarrow
X\longrightarrow 0,$$ where $h$ is a minimal right add
$(U_{t-1}\oplus P)$-approximation of $X$. Then we get that $K'\in
{\rm add}\ (\underset{0\leq j\leq t-1}{\oplus}U_j\oplus P)$ and
hence $K'\in$add $M$. So the induced sequence
$$0\longrightarrow {\rm Hom}_{A^{(m)}}(L, K')\longrightarrow {\rm Hom}_{A^{(m)}}(L, M_1')
\longrightarrow {\rm Hom}_{A^{(m)}}(L, X)\longrightarrow 0$$ is
exact for all $L\in {\rm add}\ M$. This completes the proof.
\hfill$\Box$

\vskip 0.2 in

{\bf Example 3.4.}\ We now give an example with the case of $m=1$.
Let $Q$ be the quiver $1\leftleftarrows 2$, and $Q^{(1)}$ be the
quiver $1\leftleftarrows 2\leftleftarrows 1'\leftleftarrows 2'$. Let
$A^{(1)}= kQ^{(1)}/I$ be the duplicated algebra of $kQ$. The
indecomposable projective-injective $A^{(1)}$-modules are
$P_1'= \begin{array}{c}1'\\[-2ex]22\\[-2ex]1\end{array}$ and
$P_2'= \begin{array}{c}2'\\[-2ex]11\\[-2ex]2\end{array}$
which are represented by their Loewy series.

\vskip 0.1 in

We take a generator-cogenerator $M$ of mod $A^{(1)}$ as the
following,

$$
M=1\oplus \begin{array}{c}2\\[-2ex]1\ 1 \end{array}\oplus 2'\oplus
\begin{array}{c}2'\ 2' \\[-2ex]1' \end{array}\oplus
\begin{array}{c}1'\\[-2ex]2\ 2 \\[-2ex]1 \end{array}
\oplus \begin{array}{c}2'\\[-2ex]1'\ 1' \\[-2ex]2 \end{array}
\oplus
\begin{array}{c}1'\\[-2ex]2\ 2
\end{array}\oplus \begin{array}{c}1'\ 1'\\[-2ex]2\ 2\ 2 \end{array}\oplus
\begin{array}{c}2'\ 2' \\[-2ex]1'\ 1' \ 1'\end{array} \oplus
\begin{array}{c}2'\ 2'\ 2' \\[-2ex]1'\ 1' \ 1'\ 1'\end{array}.
$$

According to Theorem 3.3, we have that ${\rm gl.dim End}_{A^{(1)}}M=
3$, and that {\rm rep.dim} $A^{(1)}=3$.

\vskip 0.1 in

If we take another generator-cogenerator $M_0$ of mod $A^{(1)}$ as
the following,

$$
M_0=1\oplus \begin{array}{c}2\\[-2ex]1\ 1 \end{array}\oplus 2'\oplus
\begin{array}{c}2'\ 2' \\[-2ex]1' \end{array}\oplus
\begin{array}{c}1'\\[-2ex]2\ 2 \\[-2ex]1 \end{array}
\oplus \begin{array}{c}2'\\[-2ex]1'\ 1' \\[-2ex]2 \end{array}
,
$$
then we have that ${\rm gl.dim End}_{A^{(1)}}M_0= 5$.

\vskip 0.2in

{\bf Remark.} \ The above example shows that the global dimension of
endomorphism algebra of a generator-cogenerator can be really bigger
than the representation dimension.

\vskip 0.2in

{\bf Theorem 3.5}\ {\it Let $A^{(m)}$ be the $m$-replicated algebra
of a hereditary algebra of $A$.  Then the dominant dimension of
$A^{(m)}$ is at least $m$. }

\vskip 0.1in

 {\bf Proof.}\ Assume that gl.dim $A^{(m)}=t.$ Then
$m+1\leq t\leq 2m+1.$ For each indecomposable projective $A$-module
$P$, we take its minimal injective resolution:
$$0\longrightarrow P\longrightarrow I_1\longrightarrow I_2
\longrightarrow \cdots \longrightarrow I_{t-2}\longrightarrow I_{t-1}
\longrightarrow\cdots \longrightarrow 0. $$
By the AR-quiver of $A^{(m)}$, we have $P,\ \Omega^{-1}_{A^{(m)}}P,\
\cdots,\  \Omega^{-(t-2)}_{A^{(m)}}P$ are not in ind $A_m$ and thus
$I_1,\ I_2,\cdots ,I_{t-1}$ are all projective-injective $A^{(m)}$.
Therefore we have a long exact sequence
$$0\longrightarrow A^{(m)}\longrightarrow N_1\longrightarrow \cdots
\longrightarrow N_{t-1} \longrightarrow \cdots \longrightarrow 0$$
with $ N_1,\cdots, N_{t-1}$ projective-injective, which implies
dom.dim $A^{(m)}\geq t-1\geq m.$ This finishes the proof.
\hfill$\Box$

\begin{description}

\item{[A]}\ M.Auslander, Representation Dimension of Artin
Algebras. Queen Mary College Mathematics Notes, Queen Mary College,
London, 1971.

\item{[ABM]}\ I.Assem, A.Beligiannis, N.Marmaridis, Right
triangulated categories with right semi-equivalences. {\it CMS
Conference Proceedings} 24(1998), 17-37.

\item{[ABST1]}\ I.Assem, T.Br$\ddot{\rm u}$stle, R.Schiffer,
G.Todorov, Cluster categories and duplicated algebras. {\it J.
Algebra.} 305(2006), 548-561.

\item{[ABST2]}\ I.Assem, T.Br$\ddot{\rm u}$stle, R.Schiffer,
G.Todorov, $m$-cluster categories and $m$-replicated algebras. {\it
Journal of pure and applied algebra.} 212(2008), 884-901.

\item{[AI]}\ I.Assem, Y.Iwanaga, On a class of
representation-finite QF-3 algebras. {\it Tsukuba J. Math. }
11(1987), 199-210.

\item{[APT]}\ I.Assem, M.I.Platzeck, S.Trepode, On the
representation dimension of tilted and laura algebras. {\it J.
Algebra.} 296(2006), 426-439.

\item{[ARS]}\ M.Auslander, I.Reiten, S.O.Smal$\phi$, \ Representation
Theory of Artin Algebras. Cambridge Univ. Press, 1995.

\item{[AS]}\ M. Auslander and S.O.Smal$\phi$, Almost split sequences in subcategories. {\it J.
Algebra.}  69(1981), 426-454.

\item{[CP]}\ F.U.Coelho, M.I.Platzeck, On the
representation dimension of some classes of algebras. {\it J.
Algebra.} 275(2004), 615-628.

\item{[EHIS]}\ K.Erdmann, T.Holm, O.Iyama, J.Schroer, Radical embedding and
representation dimension, {\it Adv.Math.} 185(2004), 159-177.

\item{[H]}\ D.Happel, Triangulated categories in the representation
theory of finite dimensional algebras. {\it Lecture Notes series
119.} Cambridge Univ. Press,1988.

\item{[HRS]}\ D.Happel, I.Reiten, S.O.Smal$\phi$, Tilting
in abelian categories and quasitilted algebras. {\it Memoris Amer.
Math.Soc.} 575(1996).

\item{[HW]}\ D.Hughes, J.Waschb$\ddot{\rm u}$sch, Trivial extensions
of tilted algebras. {\it Proc.London. Math.Soc.} 46(1983), 347-364.

\item{[IT]}\ K.Igusa, G.Todorov, On the finitistic global dimension
conjecture for artin algebras, In:Representations of algebras and
related topics,201-204.Fields Inst. Commun., 2002,45,
Amer.math.Soc., Providence, RI, 2005.

\item{[LLZ]}\ X.Lei, H.Lv, S.Zhang, Complements to the almost
complete tilting $A^{(m)}$-modules. To appear in {\it Communications
in Algebra.}

\item{[Rin]}\ C.M.Ringel, Tame algebras and integral quadratic forms.
{\it Lecture Notes in Math. 1099.} Springer Verlag, 1984.

\item{[W]}\ T. Wakamatsu, Stable equivalence of self-injective algebras and a generalization of
tilting modules. {\it J. Algebra.} 134(1990), 289-325.

\item{[X]}\ C.Xi, Representation dimension
and quasi-hereditary algebras. {\it Adv.Math.} 168(2002), 193-212.

\item{[ZZ]}\ A.Zhang, S.Zhang,  On the finitistic dimension conjecture of Artin
algebras. {\it J. Algebra.},  320(2008), 253-258.

\item{[Z1]}\ S.Zhang, Tilting mutation and duplicated
algebras. To appear in {\it Communications in Algebra.}

\item{[Z2]}\ S.Zhang, Partial tilting modules over $m$-replicated
algebras. Preprint,2008.

\end{description}

\end{document}